\documentclass[11pt,a4,fleqn]{article}
\usepackage{graphicx}
\usepackage{amsmath,amssymb,latexsym,graphics,epsfig}
\usepackage{hyperref}
\usepackage{color}
\usepackage{amsthm}

\newtheorem{theorem}{\bf Theorem}[section]
\newtheorem{proposition}[theorem]{\bf Proposition}
\newtheorem{lemma}[theorem]{\bf Lemma}
\newtheorem{corollary}[theorem]{\bf Corollary}

\newtheorem{question}{\bf Question}
\newtheorem{claim}{\bf Claim}

\newcommand{\G}{\Gamma}

\newcommand{\pf}{\noindent{\em Proof: }}
\newcommand{\epf}{\hfill\hbox{\rule{3pt}{6pt}}\\}

\numberwithin{equation}{section}

\begin{document}
\title{{\Large Variations of the spectral excess theorem for normal digraphs}}

\author{G.R. Omidi  \\[2pt]
{\small  Department of Mathematical Sciences, Isfahan University of Technology},\\
{\small Isfahan, 84156-83111, Iran}\\
{\small School of Mathematics, Institute for Research in Fundamental Sciences (IPM),}\\
{\small P.O. Box 19395-5746, Tehran, Iran }\\[2pt]
{romidi@cc.iut.ac.ir}}

\date{}

\maketitle \footnotetext[1] {This research is partially
carried out in the IPM-Isfahan Branch and in part supported
by a grant from IPM (No. 91050416).} \vspace*{-0.5cm}

\begin{abstract} It is known that every distance-regular digraph is connected and normal. An interesting question is: when is a given connected normal digraph distance-regular? Motivated by this question first we give some characterizations of weakly distance-regular digraphs.
Specially we show that whether a given connected digraph to be weakly distance-regular only depends on the equality for two invariants. Then we show that a connected normal digraph $\G$ with $d+1$ distinct eigenvalues is distance-regular  if and
only if the simple excess (the ratio of the square of mean of the numbers of shortest paths between vertices at distance $d$ to the mean of the numbers of vertices at distance $d$ from every vertex, which is zero if $d$ is greater than the diameter) is equal to the spectral excess (a number which can be computed from the spectrum of $\G$). In fact, this result is a new variation (a simple variation)
of the spectral excess theorem due to Fiol and Garigga for connected normal digraphs. Using these results we derive another variation (a weighted variation) of the spectral excess theorem for connected normal digraphs. Distance regularity of a digraph (also a graph) is in general not determined by its spectrum. For application of the simple variation we show that distance regularity of a connected normal digraph $\G$ (with $d+1$ distinct eigenvalues) is determined by its spectrum and the invariant $\delta_d$
(the mean of the numbers of vertices at distance $d$ from every vertex). Finally as an application of the weighted variation we show that every connected normal digraph $\G$ with $d+1$ distinct eigenvalues and diameter $D$ is either a bipartite digraph, or a generalized odd graph or it has odd-girth at most $\min\{2d-1,2D+1\}$, generalizing a result of van Dam and Haemers and also a recent result of Lee and Weng.

\noindent \\
\noindent{\small Keywords: The spectral excess theorem, Weakly distance-regular digraphs, Distance-regular digraphs, Normal digraphs.}\\
{\small AMS subject classification: 05C50, 05E30}

\end{abstract}

\section{Introduction}

A {\it digraph} (or a {\it directed graph}) is an ordered pair  $\G=(V,E)$ where $V$ is a set whose elements are called vertices or nodes, and
$E$ is a set of ordered pairs of vertices, called arcs or edges. In contrast, a graph where the edges are bidirectional is called an {\it undirected graph}.
A digraph having no multiple edges or loops (corresponding to a binary adjacency matrix with 0s on the diagonal) is called {\it simple}.  We will denote by $\partial(x,y)$ the distance from a vertex $x$ to a vertex $y$. For every vertex $x$ we define the   {\it directed shell}  $\Gamma^{+}_k (x)$ (resp. $\Gamma^{-}_k (x)$) the set of vertices at distance $k$ from $x$ (resp. the set of vertices from which $x$ is at distance $k$). The maximum (directed) distance between distinct pairs of vertices is called the {\it diameter} of $\G$ and is denoted by $D$. The {\it girth $g$} (resp. {\it odd-girth $g_o$}) is the smallest length of a cycle (an odd cycle) in $\G$.
In this paper, by a walk, path or cycle, we mean a directed walk, path or cycle. A digraph is {\it (strongly) connected} if there is a path between every pair of vertices. For distinct vertices $u$ and $v$ of $\G$ we say that $u$ is {\it adjacent} to $v$ if there is an edge (directed edge) from $u$ to $v$. To see more about digraphs we refer the reader to \cite{BJG}.
Throughout this paper, let $\Gamma=(V,E)$ be a connected simple digraph of order $n$ and diameter $D$.\\


The {\it adjacency matrix} $\textbf{A}=(a_{ij})$ of $\Gamma$ is the $n \times n$ matrix indexed by the vertices of $\Gamma$, with entries $a_{ij}=1$ if $i$ is adjacent to $j$, and $a_{ij}= 0$, otherwise. Digraphs with symmetric adjacency matrices are precisely the graphs. The {\it adjacency algebra} of $\Gamma$ is defined by
$\mathcal{A}(\Gamma)=\{ P(\textbf{A}): P \in C[x]\}$. A digraph $\Gamma$ is called {\it normal} if $\textbf{A}$ is a normal matrix; that is, $\textbf{A}\textbf{A}^t=\textbf{A}^t\textbf{A}$. A digraph is said to be {\it regular of degree} $k$ if $\textbf{A}J=J\textbf{A}=kJ$.
The spectrum of a digraph $\G$ is denoted by the multi-set $$Spec(\G)=\{[\lambda_0]^{m_{0}},[\lambda_1]^{m_{1}},\ldots,[\lambda_d]^{m_{d}}\}$$ where the superscripts $m_i$ denote the multiplicities of the distinct eigenvalues $\lambda_i, i = 0,1,\ldots,d$. Since $\textbf{A}$ may be not symmetric, it is possible to have non-real eigenvalues. It is known that the eigenvalue with maximum absolute value is real and simple. Here it is denoted by $\lambda_0$ and so  $m_{0}=1$. If $\Gamma$ is $k$-regular, then $\lambda_0=k$. The degree of the minimal polynomial of $\textbf{A}$ is denoted by $m_{\G}=\hat{D}+1$. Clearly $D\leq \hat{D}$, since $I, \textbf{A}, \textbf{A}^2,\ldots, \textbf{A}^D$ are linearly independent. Also $d\leq \hat{D}$ and we have the equality if and only if $\textbf{A}$ is diagonalizable (e.g., $\textbf{A}$ is symmetric or equivalently $\G$ is a graph). For $0\leq k \leq D$ the distance matrix $\textbf{A}_k$ of $\Gamma$ is defined by
$(\textbf{A}_k)_{xy}=1$ if $\partial(x,y)=k$ and $(\textbf{A}_k)_{xy}=0$, otherwise. In particular, $\textbf{A}_0=I$ and $\textbf{A}_1=\textbf{A}$.\\

For a vertex $v$, let $\delta_k(v)=|\G^{+}_k(v)|$ and $\delta'_k(v)$ be the number of paths of length $k$ from $v$ to $\G^{+}_k(v)$. Now let $\delta_k=\frac{\sum_{v\in V(\G)\delta_k(v)}}{n}$ and $\delta'_k=\frac{\sum_{v\in V(\G)\delta'_k(v)}}{n}$. A digraph is {\it geodetic} if the shortest path between any two vertices is unique. Clearly for all $0\leq k\leq D$ we have $\delta_k\leq\delta'_k$ and we have the equalities if and only if $\G$ is geodetic.  The number $\varepsilon_{\G}=\frac{\delta'^2_d}{\delta_d}$ is called the  {\it simple excess} of $\G$ (it is zero if $d>D$).
The unique polynomial $H(x)=n\frac{S(x)}{S(\lambda_0)}$ where $(x-\lambda_0)S(x)$ is the minimal polynomial of $\Gamma$ is called the {\it Hoffman polynomial}.  It is shown that, $\G$ is regular if and only if $H(\textbf{A})=J$, the all 1's matrix; see \cite{CDS}. For a vertex $v$, let $\tilde{\delta}_k(v)=\sum_{u\in\G^{+}_k(v)}(H(\textbf{A}))^2_{vu}$ and $\tilde{\delta}'_k(v)=\sum_{u\in\G^{+}_k(v)}(H(\textbf{A}))_{vu}\rho_k(v,u)$ where $\rho_k(v,u)$ is the number of paths of length $k$ from $v$ to $u$. Now let $\tilde{\delta}_k=\frac{\sum_{v\in V(\G)\tilde{\delta}_k(v)}}{n}$ and $\tilde{\delta}'_k=\frac{\sum_{v\in V(\G)\tilde{\delta}'_k(v)}}{n}$. The number $\tilde{\varepsilon}_{\G}=\frac{\tilde{\delta}'^2_d}{\tilde{\delta}_d}$  is called the  {\it weighted excess} of $\G$ (it is zero if $d>D$). Clearly for a regular digraph $\G$ we have $\delta_k(v)=\tilde{\delta}_k(v)$, $\delta'_k(v)=\tilde{\delta}'_k(v)$, $\delta_k=\tilde{\delta}_k$, $\delta'_k=\tilde{\delta}'_k$ and $\varepsilon_{\G}=\tilde{\varepsilon}_{\G}$.\\

A {\it distance-regular} graph is a regular graph such that for any two vertices $v$ and $w$ at distance $i$ the number of vertices adjacent to $w$ and at distance $j$ from $v$ is only depends on $i$ and $j$. For a connected graph with $d + 1$ distinct eigenvalues, the {\it excess} of a vertex $v$ is the number of vertices at distance $d$ from $v$. The {\it spectral excess}
of a graph is a well defined number that depends only on the spectrum
of the graph. By the spectral excess theorem we mean the remarkable result by Fiol and Garriga
\cite{FG} that a connected regular graph with $d + 1$ distinct eigenvalues is distance-regular  if and
only if the average excess is equal to the spectral excess; see \cite{v,FGG} for short proofs. This theorem is an
important tool for proving distance-regularity, and characterizing distance-regular
graphs. Since the paper \cite{FG}
appeared, different approaches (local or global) of the spectral excess theorem have been given. The interest
of the inequalities so obtained is the characterization of some kind of distance-regularity,
and this happens when equalities are attained; see \cite{DvFGG}. For more background on graph
spectra and different concepts of distance-regularity in graphs see, for instance, \cite{BCN,BH,CDS,F}.\\

The concept of a {\it distance-regular digraph} was introduced by Damerell \cite{D} in the
late 1970s. A digraph $\G$ with
diameter $D$ is  distance-regular if, for any pair of vertices $u, v$ with $\partial(u,v)=k$ for
$1\leq k\leq D$, the numbers $|\G^{+}_i(u)\cap\G^{+}_1(v)|$
for each $i$ such that $0\leq i\leq k+1$; do not depend on the chosen vertices $u$ and $v$, but
only on their distance $k$. Note that distance-regular digraphs with $g=2$ are precisely the distance-regular graphs. To see more things about distance-regular digraphs we refer the reader to \cite{EM,LM} and the references therein. If we change $\G^{+}_1(v)$ into
$\G^{-}_1(v)$ in the definition of distance-regularity, we get a new family of digraphs named by  {\it 'weakly
distance-regular digraphs'}, which are closely related to the
properties enjoyed by the distance-regular digraphs. In fact, the diameter two weakly
distance-regular digraphs are the same as the 'directed strongly regular graphs'
introduced by Duval \cite{Du}, which have recently been studied by a number of authors
and found several constructions; see \cite{Br}. The concept of weakly distance-regular digraphs was introduced by F. Comellas et al. \cite{CFGM} as a generalization of distance-regular digraphs. In fact, distance-regular digraphs are normal weakly distance-regular digraphs. Also, in \cite{CFGM} it has been shown that a digraph $\Gamma$ of diameter $D$ is weakly distance-regular if, for each nonnegative integer  $\ell\leq D$, the number $a^\ell _{uv}$ of walks of length $\ell$ from a vertex $u$ to a vertex $v$ only depends on their distance $\partial(u,v)$.\\

It is known that every distance-regular digraph is connected and normal. Here we investigate to the following interesting question:
\begin{question}\label{question1}
When is a given connected normal digraph distance-regular?
\end{question}
It has been shown that a connected digraph is distance-regular if and only if it is normal and weakly distance-regular; see \cite{CFGM}. So Question \ref{question1} is related to the following question:
\begin{question}\label{question2}
When is a given connected digraph weakly distance-regular?
\end{question}
An answer to Question \ref{question2} is given in \cite{CFGM}; see Theorem \ref{weakly distance-regular}. Motivated by this question we give further characterizations of weakly distance-regular digraphs in Section 3. In particular, we show that whether a given connected digraph to be weakly distance-regular only depends on the equality for two numbers, the first number only depends on the distance matrices and the pre-distance polynomials and the second number is the number of vertices; see Theorem \ref{wdrdg7}. To answer Question \ref{question1} we derive two new variations (a simple and a weighted variation) of the spectral excess theorem for connected normal digraphs; see Theorems \ref{drdg3} and \ref{drdg9}. We show that a connected normal digraph $\G$ is distance-regular  if and
only if the simple excess (also the weighted excess) is equal to the spectral excess (a number which can be computed from the spectrum of $\G$). So besides the spectrum, a simple
combinatorial property suffices for a connected normal digraph to be distance-regular.\\

Distance regularity of a digraph is in general not determined by the spectrum of the digraph. For application of the simple variation we show that a connected normal digraph $\G$ (with $d+1$ distinct eigenvalues) with the same spectrum and  $\delta_d$ to a distance regular digraph, is distance regular. Also, we give a new variation of the spectral excess theorem for
geodetic distance-regular digraphs and using it we show that any connected normal digraph cospectral with a geodetic distance-regular digraph, is a geodetic distance-regular digraph.\\

A distance-regular graph  with diameter $D$ and odd-girth $2D+1$ is called a {\it generalized odd graph}, also known as an {\it almost-bipartite
distance-regular graph} or a {\it regular thin near $(2D+1)$-gon}. In \cite{vH}, van Dam and Haemers showed that any connected regular graph
with $d+1$ distinct eigenvalues and odd-girth $2d+1$ is a generalized odd graph.
In a recent paper, Lee and Weng \cite{LW} used a variation of the spectral excess theorem for non-regular
graphs to show that, indeed, the regularity condition is not necessary and van Dam and Fiol \cite{vF}
gave a more direct short proof of the same result. Here as an application of the weighted variation of the spectral excess theorem for connected normal digraphs we show that every connected normal digraph $\G$ with $d+1$ distinct eigenvalues is either a bipartite digraph, or a generalized odd graph or it has odd-girth at most $\min\{2d-1,2D+1\}$. As a corollary we conclude that every connected normal digraph $\G$ with $d+1$ distinct eigenvalues and odd-girth at least $2d+1$ is a generalized odd graph, generalizing the previous mentioned result of Lee and Wang in \cite{LW}. Because the odd-girth follows from the spectrum, the latter is also a generalization of the spectral characterization of generalized odd graphs; see \cite{Huang,HLiu}.

\section{Scalar product polynomial space}
In this section, we review the concept of orthogonal polynomial related to $\G$. For more background on this concept we refer the reader to \cite{CFFG}. Let $\varphi$=$\langle m(x)\rangle$
be the ideal generated by the minimal polynomial of $\G$.
Consider the $(\hat{D}+1)$-dimensional vector space $C_{\hat{D}}[x]\simeq C[x]/\varphi$  with the following scalar product.
\begin{equation}\label{eq(2)}
\langle f , g\rangle = \dfrac{1}{n} tr( f(\textbf{A})g(\textbf{A})^\ast)
\end{equation}
where $g(\textbf{A})^\ast$ is the transpose of the conjugate of $g(\textbf{A})$
and norm defined by $\parallel f \parallel = \sqrt{\langle f,f \rangle}$.
It is clear that  $\parallel f \parallel=0$ if and only if $f=0$. The {\it projection} of $f(x)$ into $g(x)$ is defined by
$$Proj_{g}(f)=\dfrac{\langle g,f\rangle}{\langle g,g\rangle} g(x)$$
where $\langle g,f\rangle$ denotes the scalar product in \ref{eq(2)}. The projection of $f$ into $\{g_k\}_{k=0}^{l}$ is defined by
$\sum_{k=0}^{l} Proj_{g_k}(f).$
Clearly $\lbrace 1,x,x^2,\dots x^{\hat{D}} \rbrace$ is a basis of $C_{\hat{D}}[x]$. We use the Gram-Schmidt process and generate an orthogonal set
$\lbrace \mathcal{P}_0, \mathcal{P}_1,\ldots, \mathcal{P}_{\hat{D}} \rbrace$ that spans the same $(\hat{D}+1)$-dimensional space as follows.
Let $\mathcal{P}_0 = 1$, $\mathcal{P}_1= x - \dfrac{\langle 1,x \rangle}{\langle1,1\rangle} 1 = x$ and for $0\leq i\leq \hat{D}-1$
$$\mathcal{P}_{i+1}= x^{i+1}-\sum_{k=0}^i Proj _{\mathcal{P}_{k}}(x^{i+1}).$$
One can easily see that the polynomials $\{\mathcal{P}_0, \mathcal{P}_1,\ldots, \mathcal{P}_{\hat{D}} \}$ is an orthogonal basis of $C_{\hat{D}}[x]$ and
$\mathcal{P}_i$ is monic (the leading coefficient  is equal to $1$) and it has degree $i$. These polynomials are referred to as the
{\it pre-distance polynomials} of $\G$. \\

For two $n\times n$ matrices $\textbf{C}$ and $\textbf{D}$ assume that $\langle \textbf{C} , \textbf{D}\rangle = \dfrac{1}{n} tr( \textbf{C}\textbf{D}^\ast)$ and $$Proj_{\textbf{D}}(\textbf{C})=\dfrac{\langle \textbf{C},\textbf{D}\rangle}{\langle \textbf{D},\textbf{D}\rangle} \textbf{D}.$$ The projection of $\textbf{C}$ into a set of matrices is the summation of the projections of $\textbf{C}$ into all matrices in that set. From the definitions for every $f,g\in C_{\hat{D}}[x]$ we have $\langle f , g\rangle=\langle f(\textbf{A}) , g(\textbf{A})\rangle$, $\delta_k=\langle \textbf{A}_k , \textbf{A}_k\rangle$ and $\varepsilon_{\G}=\frac{\langle \textbf{A}_d , \textbf{A}^d\rangle^2}{\delta_d}$ if $d\leq D$ (it is zero if $d> D$).
Now assume that the weighted distance matrix $\tilde{\textbf{A}}_k$ is an $n\times n$ matrix with $(\tilde{\textbf{A}}_k)_{uv}=(H(\textbf{A}))_{uv}(\textbf{A}_k)_{uv}$. Clearly $H(\textbf{A})=\tilde{\textbf{A}}_0+\tilde{\textbf{A}}_1+\cdots+\tilde{\textbf{A}}_D$. From the definitions we have $\tilde{\delta}_k=\langle \tilde{\textbf{A}}_k , \tilde{\textbf{A}}_k\rangle$ and $\tilde{\varepsilon}_{\G}=\frac{\langle \tilde{\textbf{A}}_d , \textbf{A}^d\rangle^2}{\tilde{\delta}_d}$ if $d\leq D$ (it is zero if $d> D$). In the sequel, for every $0\leq k\leq \hat{D}$ we assume that $\mathcal{Q}_k=\sum_{i=0}^{k} \mathcal{P}_i$ and by $\{P_0, P_1,\ldots, P_{\hat{D}}\}$ we mean the set of polynomials obtained by normalizing the pre-distance polynomials so that $\|P_k\|^2=\delta_k$ for $0\leq k\leq D$ and $P_k=\mathcal{P}_k$ for $D< k\leq \hat{D}$. For a weakly distance- regular digraph $\Gamma$ we have
$P_i(\textbf{A})=\textbf{A}_i$ for $0\leq i \leq \hat{D}=D$, $$\langle P_i,P_j\rangle=\dfrac{1}{n} tr( \textbf{A}_{i} \textbf{A}^T_{j})= \dfrac{1}{n} \sum_{u} \vert \Gamma^{+}_i (u)\cap \Gamma^{+}_j (u) \vert =0$$ for $i\neq j$ and $\|P_{i}\|^2=\dfrac{1}{n} tr(\textbf{A}_{i} \textbf{A}^T_{j})=\delta_k= P_k(\Delta).$\\

Through this paper, we assume that  $\G$ is a connected simple digraph on $n$ vertices with $d+1$ distinct eigenvalues, diameter $D$, $m_{\G}=\hat{D}+1$, distance matrices $\{\textbf{A}_k\}_{k=0}^{D}$, weighted distance matrices $\{\tilde{\textbf{A}}_k\}_{k=0}^{D}$, pre-distance polynomials $\{\mathcal{P}_k\}_{k=0}^{\hat{D}}$ and normalized pre-distance polynomials $\{P_k\}_{k=0}^{\hat{D}}$.

\section{Weakly distance-regular digraphs}

In this section, we give some characterizations of weakly distance-regular digraphs. Then we conclude that whether a given connected digraph to be weakly distance-regular only depends on the equality for two invariants. The following theorem gives a characterization for weakly distance-regular digraphs.

\begin{theorem}\label{weakly distance-regular}{\cite{CFGM}}\\
For a connected digraph $\Gamma$ the following are equivalent:
\begin{itemize}\item[\rm(i)]
$\G$ is a weakly distance-regular digraph;
\item[\rm(ii)]  The distance matrix $\textbf{A}_k$ is a polynomial of degree $k$ in the adjacency matrix $\textbf{A}$; that
is, $\textbf{A}_k=P_k(\textbf{A})$ for each $k= 0,1,\ldots, D$, where $P_k\in Q[x]$;
\item[\rm(iii)]  The set of distance matrices $\lbrace \textbf{A}_k: k= 0,1,\ldots, D \rbrace$
 is a basis of the adjacency algebra $\mathcal{A}(\G)$;
 \item[\rm(iv)]  For any two vertices $u, v$  at distance $\partial (u,v)=k$; the numbers
$$p^k_{ij}(u,v)=|\G^{+}_i(u)\cap \G^{-}_j(v)|$$
do not depend on the vertices $u$ and $v$, but only on their distance $k$; in which case
they are denoted by $p^k_{ij}$ (note that $p^k_{ij}=0$ when $k>i+j$).
\end{itemize}
\end{theorem}

\begin{lemma}\label{wdrdg1} For a connected digraph $\Gamma$ we have
\begin{equation}\label{equ1}
\sum_{k=0}^{D}\frac{\langle \textbf{A}_k,P_k(\textbf{A}) \rangle^2}{\delta_k}\leq n,
\end{equation}
and the equality holds if and only if $\G$ is weakly distance-regular.
\end{lemma}

\pf Set $\eta_k=\delta_k-\frac{\langle \textbf{A}_k,P_k(\textbf{A}) \rangle^2}{\delta_k}$ for $0\leq k\leq D$. The projection of $P_k(\textbf{A})$ into $\textbf{A}_k$ is:
$$\bar{P}_k(\textbf{A})=\frac{\langle \textbf{A}_k,P_k(\textbf{A}) \rangle}{\| \textbf{A}_k\|^2}\textbf{A}_k.$$
Now consider the equality $P_k(\textbf{A})=\bar{P}_k(\textbf{A})+N_k$, where $\langle N_k,\bar{P}_k(\textbf{A})\rangle=0$. Then, from Pythagores theorem and the above equality, we have

\begin{equation}\label{equ2}
\| N_k\|^2=\| P_k(\textbf{A})\|^2-\| \bar{P}_k(\textbf{A}) \|^2=\delta_k-\frac{\langle \textbf{A}_k,P_k(\textbf{A}) \rangle^2}{\delta_k}\geq 0.
\end{equation}

Therefore $\eta_k\geq 0$ and so $\sum_{k=0}^{D}\eta_k\geq 0$. Hence
$$\sum_{k=0}^{D}\frac{\langle \textbf{A}_k,P_k(\textbf{A}) \rangle^2}{\delta_k}\leq \sum_{k=0}^{D}\delta_k=n.$$

To give a proof for the second part, first assume that $\G$ is weakly distance-regular. Then $\textbf{A}_k=P_k(\textbf{A})$ for $0\leq k\leq D$ and so $\langle \textbf{A}_k,P_k(\textbf{A}) \rangle=\delta_k$. Hence $\eta_k=0$ for each $k$. Therefore $$n=\sum_{k=0}^{D}\frac{\langle \textbf{A}_k,P_k(\textbf{A}) \rangle^2}{\delta_k}.$$

Now assume that the equality in \ref{equ1} is attained. Then $\sum_{k=0}^{D}\eta_k=0$ and so using \ref{equ2}, $\eta_k= 0$ for each $k$. Clearly $\eta_k=0$ if and only if $N_k=0$ and hence $\bar{P}_k(\textbf{A})=P_k(\textbf{A})=\frac{\langle \textbf{A}_k,P_k(\textbf{A}) \rangle}{\delta_k}\textbf{A}_{k}$. Since $\textbf{A}_k$ is a
polynomial of degree $k$ in $\textbf{A}$, using Theorem \ref{weakly distance-regular} $\G$ is a weakly distance-regular digraph.\epf

The following corollary is an immediate consequent from the proof of Lemma \ref{wdrdg1}.

\begin{corollary}\label{wdrdg2} For every  $0\leq k\leq D$ we have, $$\langle \textbf{A}_k,P_k(\textbf{A}) \rangle\leq \delta_k.$$
All equalities hold if and only if $\G$ is weakly distance-regular.
\end{corollary}

 By considering the projections of $\mathcal{P}_k(\textbf{A})$ into $\textbf{A}_k$ for all  $0\leq k\leq D$ and applying  the same argument in Lemma \ref{wdrdg1} we get the following characterization of weakly distance-regular digraphs which are geodetic.

\begin{theorem}\label{gwdrdg1} For a connected digraph $\G$ we have $n\leq\| \mathcal{Q}_D\|^2,$
and the equality holds if and only if $\G$ is a geodetic weakly distance-regular digraph.
\end{theorem}

\begin{lemma}\label{wdrdg5} For a connected digraph $\Gamma$ we have
\begin{equation}\label{equ4}
\sum_{k=0}^{D}\sum_{j=k}^{D}\frac{\langle \textbf{A}_k,P_j(\textbf{A}) \rangle^2}{\delta_j}\leq n
\end{equation}
and the equality holds if and only if $\G$ is weakly distance-regular.
\end{lemma}

\pf Set $\eta_k=\delta_k-\sum_{j=k}^{D}\frac{\langle \textbf{A}_k,P_j(\textbf{A}) \rangle^2}{\delta_j}$. The projection of $\textbf{A}_k$ into $\{P_j(\textbf{A})\}_{j=0}^{D}$ is:
$$\bar{\textbf{A}}_k=\sum_{j=0}^{D}\frac{\langle \textbf{A}_k,P_j(\textbf{A}) \rangle}{\delta_j}P_j(\textbf{A}) =\sum_{j=k}^{D}\frac{\langle \textbf{A}_k,P_j(\textbf{A}) \rangle}{\delta_j}P_j(\textbf{A}).$$
Now consider the equality $\textbf{A}_k=\bar{\textbf{A}}_k+N_k$, where $\langle N_k,\bar{\textbf{A}}_k\rangle=0$. Then, from Pythagores theorem and the above equality, we have

\begin{equation}\label{equ5}
\| N_k\|^2=\| \textbf{A}_k\|^2-\| \bar{\textbf{A}}_k \|^2=\delta_k-\sum_{j=k}^{D}\frac{\langle \textbf{A}_k,P_j(\textbf{A}) \rangle^2}{\delta_j}=\eta_k\geq 0.
\end{equation}

Therefore $\sum_{k=0}^{D}\eta_k\geq 0$ and hence
$$\sum_{k=0}^{D}\sum_{j=k}^{D}\frac{\langle \textbf{A}_k,P_j(\textbf{A}) \rangle^2}{\delta_j}\leq \sum_{k=0}^{D}\delta_k=n.$$

If $\G$ is weakly distance-regular, then $\textbf{A}_i=P_i(\textbf{A})$ for $0\leq i\leq D$ and so the equality in \ref{equ4} is attained.
Now assume that the equality in \ref{equ4} is attained. Then $\sum_{k=0}^{D}\eta_k=0$ and so using \ref{equ5}, $\eta_k= 0$ for each $k$.
Clearly $\eta_k=0$ if and only if $\textbf{A}_k=\sum_{j=k}^{D}\frac{\langle \textbf{A}_k,P_j(\textbf{A}) \rangle}{\delta_j}P_j(\textbf{A}).$ Therefore,
$$\begin{bmatrix}
 \textbf{A}_1 \\
\textbf{A}_2\\
\vdots \\
\textbf{A}_D
\end{bmatrix} =\begin{bmatrix}
 a_{11} & a_{12} & \ldots & a_{1(D-1)} &  a_{1D} \\
0 & a_{22} &  \ldots & a_{2(D-1)} &  a_{2D}\\\vdots & \vdots & \ldots & \vdots   & \vdots\\
0 & 0 & \ldots & 0 & a_{DD}
\end{bmatrix}.\begin{bmatrix}
 P_1(\textbf{A}) \\
P_2(\textbf{A})\\\vdots \\
P_D(\textbf{A})
\end{bmatrix}
$$
where $\textbf{C}=[a_{ij}]$ for $a_{ij}=\frac{\langle \textbf{A}_i,P_j(\textbf{A}) \rangle}{\delta_j}$
 is an upper triangular matrix. Since $\langle \textbf{A}_i,\textbf{A}^i \rangle >0$ for each $0\leq i\leq D$, all
$a_{ii}=\frac{\langle \textbf{A}_i,P_i(\textbf{A}) \rangle}{\delta_i}>0$ for any $0\leq i\leq D$; the matrix
$C$ is non-singular and its inverse $C^{-1}=[b_{ij}]$ is also an upper  triangular matrix. Therefore, $P_k(\textbf{A})=\sum_{j=k}^{D}b_{kj}\textbf{A}_j.$ On the other hand, $(P_k(\textbf{A}))_{uv}=0$ if $\partial(u,v)>k$ and so $b_{kj}=0$ for $j>k$. Hence $\textbf{A}_k$ is a
polynomial of degree $k$ in $\textbf{A}$ and so using Theorem \ref{weakly distance-regular} $\G$ is a weakly distance-regular digraph.\epf

\begin{corollary}\label{wdrdg6} For every  $0\leq k\leq D$, $$\sum_{j=k}^{D}\frac{\langle \textbf{A}_k,P_j(\textbf{A}) \rangle^2}{\delta_j}\leq \delta_k$$
and all the equalities hold if and only if $\G$ is weakly distance-regular.
\end{corollary}

Here we generalize Lemmas \ref{wdrdg1} and \ref{wdrdg5} as follow.
\begin{theorem}\label{wdrdg7} For every $0\leq k\leq D$ assume that $S_k\subseteq \{0,1,\ldots,D\}$ is non-empty. Then
\begin{itemize}\item[\rm(i)]
\begin{equation}\label{equ7}
\sum_{k=0}^{D}\sum_{j\in S_k}\frac{\langle \textbf{A}_j,P_k(\textbf{A}) \rangle^2}{\delta_j}\leq n.
\end{equation}
\item[\rm(ii)] If $k\in S_k$ for every $0\leq k\leq D$, then
\begin{equation}\label{equ6}
\sum_{k=0}^{D}\sum_{j\in S_k}\frac{\langle \textbf{A}_k,P_j(\textbf{A}) \rangle^2}{\delta_j}\leq n.
\end{equation}
In each case the equality holds if and only if $\G$ is weakly distance-regular.
\end{itemize}
\end{theorem}
\pf
In the proof of Lemma \ref{wdrdg1} if we consider the projection of $P_k(\textbf{A})$ into $\{\textbf{A}_i\}_{i\in S_k}$, $S_k\subseteq \{0,1,\ldots,D\}$, and continue the proof we get the Inequality \ref{equ7}. In the proof of Lemma \ref{wdrdg5} if we consider the projection of $\textbf{A}_k$ into $\{P_i(\textbf{A})\}_{i\in S_k}$, $k\in S_k\subseteq \{0,1,\ldots,D\}$, and continue the proof we get the Inequality \ref{equ6}.

Arguments similar to the proofs of Lemmas \ref{wdrdg1} and \ref{wdrdg5} show that in each case the equality holds if and only if $\G$ is weakly distance-regular.
\epf


The following theorem is a weighted variation of Theorem \ref{weakly distance-regular} that we will use later
on.

\begin{theorem}\label{weakly distance-regular2}
The following are equivalent:
\begin{itemize}\item[\rm(i)] $\G$ is a weakly distance-regular digraph;
\item[\rm(ii)] The matrix $\tilde{\textbf{A}}_k$ is a polynomial of degree $k$ in the adjacency matrix $\textbf{A}$; that
is, $\tilde{\textbf{A}}_k=Q_k(\textbf{A})$ for each $k= 0,1,\ldots, D$, where $Q_k\in Q[x]$;
\item[\rm(iii)] The set of matrices $\lbrace \tilde{\textbf{A}}_k: k= 0,1,\ldots, D \rbrace$
 is a basis of the adjacency algebra $\mathcal{A}(\Gamma)$;
\item[\rm(iv)] For any two vertices $u, v$  at distance $\partial (u,v)=k$; the numbers
$$\tilde{p}^k_{ij}(u,v)=\sum_{w\in \G^{+}_i(u)\cap \G^{-}_j(v)} (H(\textbf{A}))_{ww}$$
do not depend on the vertices $u$ and $v$, but only on their distance $k$; in which case
they are denoted by $\tilde{p}^k_{ij}$ (note that $\tilde{p}^k_{ij}=0$ when $k>i+j$).\end{itemize}
\end{theorem}
\pf $(i) \Longleftrightarrow (ii)$: If $\G$ is weakly distance-regular, then $\textbf{A}_k=\tilde{\textbf{A}}_k=P_k(\textbf{A})$ for each $k= 0,1,\ldots, D$. Now let $\tilde{\textbf{A}}_k=Q_k(\textbf{A})$ for each $k= 0,1,\ldots, D$. Since $(\textbf{A}-\lambda_0I)H(\textbf{A})=0$, each column of $H(\textbf{A})$ is a multiple of $\alpha$ where $\alpha$ is the Pron-eigenvector of $\textbf{A}$. So we may assume that $H(\textbf{A})=[\beta_1\alpha;\beta_2\alpha;\cdots;\beta_{n}\alpha;]$.
On the other hand, $(\textbf{A}^t-\lambda_0I)H^t(\textbf{A})=0$ and so each row of $H(\textbf{A})$ is an eigenvector of $\textbf{A}^t$ corresponding to $\lambda_0$. Therefore, $H(\textbf{A})=\alpha\beta^t$ where $\beta$ is an eigenvector of $\textbf{A}^t$ corresponding to $\lambda_0$. Note that $\textbf{A}^t$ is the adjacency matrix of a digraph obtained from $\G$ by reversing the direction of all edges and so all entries of $\beta$ have the same sign. Since $\tilde{\textbf{A}}_0=Q_{0}(\textbf{A})=d_0I$, we have $\alpha_v\beta_v=d_0$ for each vertex $v$. On the other hand, $\tilde{\textbf{A}}_1=Q_{1}(\textbf{A})=d_1\textbf{A}$ and so for each edge $(u,v)$ we have $\alpha_u\beta_v=d_1$ and so $\alpha_u=\frac{d_1}{d_0}\alpha_v$. Since $\G$ is connected each edge $e=(u,v)$ lies on a cycle $C_e=u_1u_2\cdots u_{t_e}u_1$ where $u_1=u$ and $u_{2}=v$ and so $$\alpha_{u_1}=\frac{d_1}{d_0}\alpha_{u_2}=(\frac{d_1}{d_0})^2\alpha_{u_3}=\cdots=(\frac{d_1}{d_0})^{t_e}\alpha_{u_1}.$$ Therefore $(\frac{d_1}{d_0})^{t_e}=1$, hence $\frac{d_1}{d_0}=1$ ($d_0$ and $d_1$ have the same sign) and so $\alpha_u=\alpha_v$ for each edge $e=(u,v)$. Since $\G$ is connected, all entries of $\alpha$ are equal and so $\G$ is regular. Therefore $H(\textbf{A})=J$ and $\textbf{A}_k=\tilde{\textbf{A}}_k=Q_{k}(\textbf{A})$. By Theorem \ref{weakly distance-regular} $\G$ is weakly distance-regular.

$(ii) \Longrightarrow (iii)$: Clearly $\lbrace \tilde{\textbf{A}}_k: k= 0,1,\ldots, D \rbrace$ constitutes a set of linearly independent
matrices of the adjacency algebra $\mathcal{A}(\G)$. Moreover, since  $$H(\textbf{A})=\tilde{\textbf{A}}_0+\tilde{\textbf{A}}_1+\cdots+\tilde{\textbf{A}}_D,$$
$\tilde{\textbf{A}}_k=Q_k(\textbf{A})$ and $(\textbf{A}-\lambda_0I)H(\textbf{A})=0$, it follows that the minimum polynomial
of $\G$, $(x-\lambda_0)H(x)$, has degree $m_{\G}=\hat{D}+1=D+1$ and $\dim (\mathcal{A}(\G))=D+1$. Hence, $\lbrace \tilde{\textbf{A}}_k: k= 0,1,\ldots, D \rbrace$ is a basis of $\mathcal{A}(\G)$.

$(iii) \Longrightarrow (iv)$: Assume that $u$ and $v$ are two vertices of $\G$ with $\partial(u,v)=k$.
Note that $(\tilde{\textbf{A}}_i\tilde{\textbf{A}}_j)_{uv}=\alpha_u\beta_v \tilde{p}^k_{ij}(u,v)$. On the other hand,  because of (iii), $\tilde{\textbf{A}}_i\tilde{\textbf{A}}_j$, is a linear combination of the basis $\lbrace \tilde{\textbf{A}}_k: k= 0,1,\ldots, D \rbrace$, say,
$$\tilde{\textbf{A}}_i\tilde{\textbf{A}}_j=\sum_{l=0}^{D}\gamma^l_{ij}\tilde{\textbf{A}}_l.$$
Consequently, $\alpha_u\beta_v \tilde{p}^k_{ij}(u,v)=\alpha_u\beta_v\gamma^k_{ij}$ or equivalently $\tilde{p}^k_{ij}(u,v)=\gamma^k_{ij}$ for any two vertices $u, v$ at distance $k$.

$(iv) \Longrightarrow (i)$: Since the number
$\tilde{p}^0_{00}(w,w)$ does not depend on the vertex $w$, we conclude that $a=(H(\textbf{A}))_{ww}$ is independent to the vertex $w$. Hence

$$\frac{\tilde{p}^k_{ij}(u,v)}{a}=|\G^{+}_i(u)\cap \G^{-}_j(v)|$$

does not depend on the vertices $u$ and $v$, but only on their distance $k$.  The assertion holds by  Theorem \ref{weakly distance-regular}.

\epf

\section{New variations of the spectral excess theorem}

In this section, based on some results in Section 3 we derive two new variations (a simple and a weighted variation) of the spectral excess theorem for connected normal digraphs. In fact, we show that a connected normal digraph $\G$ is distance-regular  if and
only if the simple excess (also the weighted excess) is equal to the spectral excess (a number which can be computed from the spectrum of $\G$).

Let $\textbf{A}$ be normal with eigenvalues $\{\theta_1, \theta_2,\ldots, \theta_n\}$. Then by Theorem 1.1 in \cite{CFGM} $\textbf{A}^t=\textbf{A}^\ast=f(\textbf{A})\in \mathcal{A}(\G)$ and $\textbf{A}$ diagonalize by means of a unitary matrix; that is, $\textbf{U}^\ast \textbf{A}\textbf{U}=\textbf{D}$ for some matrix
$\textbf{U}$ such that $\textbf{U}\textbf{U}^\ast=I$; and $\textbf{D}=diag(\theta_1; \theta_2;\ldots; \theta_n)$. Therefore $\hat{D}=d$ and
$$\textbf{A}^t=\textbf{A}^\ast=\textbf{U} \textbf{D}^\ast \textbf{U}^\star=f(\textbf{A})=\textbf{U} f(\textbf{D})\textbf{U}^\ast.$$
Hence $\textbf{D}^\ast=f(\textbf{D})=diag(\bar{\theta}_1; \bar{\theta}_2;\ldots; \bar{\theta}_n)=diag(f(\theta_1); f(\theta_2);\ldots; f(\theta_n))$ and so $\bar{\theta}_i=f(\theta_i)$ for each $1\leq i\leq n$ where $\bar{\theta}_i$ is the conjugate of $\theta_i$. Now let $$Spec(\G)=\{[\lambda_0]^{m_{0}},[\lambda_1]^{m_{1}},\ldots,[\lambda_d]^{m_{d}}\}$$ where the superscripts $m_i$ denote the multiplicities of the distinct eigenvalues $\lambda_i, i = 0,1,\ldots,d$. Using Lagrange's interpolation formula we have
\begin{equation}\label{eqn}
f(x)=\sum_{j=0}^{d} \bar{\lambda}_j\prod_{i\neq j} \frac{x-\lambda_i}{\lambda_j-\lambda_i}.
\end{equation}
Using Equation \ref{eqn}, all coefficients of $f(x)$ can be computed from the spectrum of $\G$. Based on this fact if all coefficients of the two polynomials $p(x)$ and $q(x)$ are real and they can also be expressed in terms of the spectrum of $\G$, then $$\langle p(x),q(x)\rangle=\dfrac{1}{n} tr( p(\textbf{A})q(f(\textbf{A})))=\dfrac{1}{n} \sum_{i=0}^{d}m_ip(\lambda_i)q(f(\lambda_i))$$ can be computed from the spectrum of $\G$.
The following proposition shows that the pre-distance polynomials of a given connected normal digraph are determined by its spectrum. The proof follows from the above facts.

\begin{proposition}\label{pre-distance} All coefficients of the pre-distance polynomials of a given connected normal digraph $\G$ can be computed from the spectrum of $\G$.
\end{proposition}
\pf
Based on the previous discussion the assertion holds by using induction and the following recursive equations for $0\leq i\leq d-1$.
$$\mathcal{P}_{i+1}(x)= x^{i+1}-\sum_{k=0}^i \dfrac{\langle x^{i+1}, \mathcal{P}_{k}(x)\rangle}{\langle \mathcal{P}_{k}(x),\mathcal{P}_{k}(x)\rangle} \mathcal{P}_{k}(x).$$
\epf

By Proposition \ref{pre-distance} for every $0\leq k\leq d$, $\varepsilon_{k}=\langle \mathcal{P}_k(x),\mathcal{P}_k(x)\rangle$ can be computed from the spectrum of $\G$.
The number  $\varepsilon_{d}$ is called {\it the spectral excess} of $\G$. The spectral excess has an important role in characterizing the distance-regular digraphs in this section. The following result is a simple variation of the spectral excess theorem for connected normal graphs.

\begin{theorem}[A simple variation of spectral excess theorem]\label{drdg3} Let $\G$ be a connected digraph. Then
$$\varepsilon_{\G}\leq \varepsilon_{d}$$
and the equality holds for a connected normal digraph $\G$ if and only if $\G$ is a distance-regular digraph.
\end{theorem}
\pf The inequality $\varepsilon_{\G}\leq \varepsilon_{d}$ follows from $\parallel \mathcal{P}_d(\textbf{A})- Proj_{\textbf{A}_d}(\mathcal{P}_d(\textbf{A}))\parallel\geq 0$ when $d\leq D$ and the fact that $\varepsilon_{\G}=0$ if $d>D$. If $\G$ is distance-regular,  then $\G$ is a connected normal digraph with $D=d$ and $\textbf{A}_D=P_D(\textbf{A})$. So $\varepsilon_{\G}=\varepsilon_{d}$.

Now assume that for a connected normal digraph $\G$ we have
$\varepsilon_{\G}=\varepsilon_{d}$. We show that $\G$ is a distance-regular digraph. Clearly $D=d$ (since otherwise $\varepsilon_{\G}=0$ and $\varepsilon_{d}>0$, a contradiction) and $\langle \textbf{A}_D,P_D(\textbf{A}) \rangle=\delta_D$. Therefore, $\parallel P_D(\textbf{A})- Proj_{\textbf{A}_D}(P_D(\textbf{A}))\parallel= 0$ and so $\textbf{A}_D=P_D(\textbf{A})$. Since $\textbf{A}$ is normal, $\textbf{A}^t$ is in the adjacency algebra and so for each $0\leq k\leq D$, $\textbf{A}^tP_k(\textbf{A})$ is in the adjacency algebra. Therefore,
$$\textbf{A}^tP_k(\textbf{A})=\sum_{i=0}^{D}\frac{\langle \textbf{A}^tP_k(\textbf{A}),P_i(\textbf{A}) \rangle}{\delta_i}P_i(\textbf{A}).$$
One can easily see that if $(P_k(\textbf{A}))_{uv}=0$ when $(\textbf{A}_k)_{uv}=0$, then the $(u,v)$-th entry of $\textbf{A}^tP_k(\textbf{A})$ is zero if $\partial(u,v)\leq k-2$ and so we have $\langle \textbf{A}^tP_k(\textbf{A}),P_i(\textbf{A}) \rangle=0$ for $i\leq k-2$. Therefore, if $\tilde{\textbf{A}}_k=c_kP_k(\textbf{A})$ for some constant $c_k$, then

\begin{equation}\label{eq8}
\textbf{A}^tP_k(\textbf{A})=\sum_{i=k-1}^{D}\frac{\langle \textbf{A}^tP_k(\textbf{A}),P_i(\textbf{A}) \rangle}{\delta_i}P_i(\textbf{A}).
\end{equation}
Since $\textbf{A}_D=P_D(\textbf{A})$, using Equation \ref{eq8} for $k=D$ we conclude that $(P_{D-1}(\textbf{A}))_{uv}$ is zero if $\partial(u,v)\neq D-1$ (note that $\langle \textbf{A}^tP_D(\textbf{A}),P_{D-1}(\textbf{A}) \rangle=\langle \textbf{A}^t\textbf{A}_D,P_{D-1}(\textbf{A}) \rangle\neq 0$). On the other hand, $H(\textbf{A})$ is in adjacency algebra. Therefore
 $$H(\textbf{A})=\tilde{\textbf{A}}_0+\tilde{\textbf{A}}_1+\cdots+\tilde{\textbf{A}}_D=\sum_{i=0}^{D}\frac{\langle H(\textbf{A}),P_i(\textbf{A}) \rangle}{\delta_i}P_i(\textbf{A}).$$
Using the fact that $\textbf{A}_D=P_D(\textbf{A})$ and the $(u,v)$-entry of $P_i(\textbf{A})$ is zero if $\partial(u,v)>i$, we have $\tilde{\textbf{A}}_{D-1}=\frac{\langle H(\textbf{A}),P_{D-1}(\textbf{A}) \rangle}{\delta_{D-1}}P_{D-1}(\textbf{A})$ and $\tilde{\textbf{A}}_{D}=\frac{\langle H(\textbf{A}),P_{D}(\textbf{A}) \rangle}{\delta_{D}}P_{D}(\textbf{A})$.

Now we show that $\tilde{\textbf{A}}_k=\frac{\langle H(\textbf{A}),P_{k}(\textbf{A}) \rangle}{\delta_{k}}P_{k}(\textbf{A})$ for each $0\leq k\leq D$ by induction on $D-k$ and hence using Theorem \ref{weakly distance-regular2}, $\G$ is a distance-regular digraph. Now let $\tilde{\textbf{A}}_i=\frac{\langle H(\textbf{A}),P_{i}(\textbf{A}) \rangle}{\delta_{i}}P_{i}(\textbf{A})$ for $r\leq i\leq D$. Using Equation \ref{eq8} for $k=r$, we conclude that the $(u,v)$-th entry of $P_{r-1}(\textbf{A})$ is zero if $\partial(u,v)\leq r-2$ and so
$(P_{r-1}(\textbf{A}))_{uv}$ is zero if $\partial(u,v)\neq r-1$ (note that $\langle \textbf{A}^tP_r(\textbf{A}),P_{r-1}(\textbf{A}) \rangle\neq 0$). Using the equation
 $$\tilde{\textbf{A}}_0+\tilde{\textbf{A}}_1+\cdots+\tilde{\textbf{A}}_{r-1}=\sum_{i=0}^{r-1}\frac{\langle H(\textbf{A}),P_i(\textbf{A}) \rangle}{\delta_i}P_i(\textbf{A})$$
and the fact that the $(u,v)$-entry of $P_i(\textbf{A})$ is zero if $\partial(u,v)>i$, we have $\tilde{\textbf{A}}_{r-1}=\frac{\langle H(\textbf{A}),P_{r-1}(\textbf{A}) \rangle}{\delta_{r-1}}P_{r-1}(\textbf{A})$, which completes the proof.\epf

Since $\varepsilon_{d}$ can be computed from the spectrum the following result is an immediate consequent from Theorem \ref{drdg3}.

\begin{theorem} Any connected normal digraph $\G$ cospectral with a distance-regular digraph $\G'$ with $\delta_d(\G')=\delta_d(\G)$, is a distance-regular digraph.
\end{theorem}
\pf
Since $\delta_d(\G')=\delta_d(\G)>0$ clearly for $\G$ we have $D=d$  and $\delta_D=\frac{\mathcal{P}^2_D(\lambda_0)}{\langle \mathcal{P}_D , \mathcal{P}_D\rangle}$
only depends on the spectrum of $\G$. Also, $\G$ is regular and so  $H(\textbf{A})=J$.
Consider the vertices $u$ and $v$ of $\G$ with $\partial (u,v)=D$. Since $H(x)$ has degree $D$, we have $(\textbf{A}^{D})_{uv}=\frac{\pi_0}{n} (H(\textbf{A}))_{uv}=\frac{\pi_0}{n}$ where $\pi_0=\prod_{i=1}^{d}(\lambda_0-\lambda_i)$ and so $\varepsilon_{\G}=\frac{\langle \textbf{A}_D , \textbf{A}^D\rangle^2}{\delta_D}=(\frac{\pi_0}{n})^2\delta_D$
only depends on the spectrum of $\G$. On the other hand, $\varepsilon_{d}$ depends also on the spectrum of $\G$. Therefore $\varepsilon_{\G}$ must equal to $\varepsilon_{d}$, because it does so for $\G'$. Hence by Theorem \ref{drdg3} $\G$ is distance regular.
\epf

Since $\mathcal{Q}_d=\sum_{k=0}^{d}\varepsilon_{k}$ follows from the spectrum using Theorem \ref{gwdrdg1} we get a new variation of spectral excess theorem for
geodetic distance-regular digraphs.

\begin{theorem}\label{gdrdg1} A connected normal digraph $\G$ is a geodetic distance-regular digraph if and only if
$\| \mathcal{Q}_d\|^2=n$.
\end{theorem}

Again since $\mathcal{Q}_d$ can be computed from the spectrum Theorem \ref{gdrdg1} implies the following corollary.

\begin{corollary} Any connected normal digraph cospectral with a geodetic distance-regular digraph, is a geodetic distance-regular digraph.
\end{corollary}

Now we present a weighted variation of spectral excess theorem for connected normal digraphs.

\begin{theorem}[A weighted variation of spectral excess theorem]\label{drdg9} Let $\G$ be a digraph. Then
$$\tilde{\varepsilon}_{\G}\leq \varepsilon_{d}$$
and the equality holds for a connected normal digraph $\G$ if and only if $\G$ is distance-regular.
\end{theorem}

\pf The inequality $\tilde{\varepsilon}_{\G}\leq \varepsilon_{d}$ follows from $\parallel \mathcal{P}_d(\textbf{A})- Proj_{\tilde{\textbf{A}}_d}(\mathcal{P}_d(\textbf{A}))\parallel\geq 0$ when $d\leq D$ and the fact that $\tilde{\varepsilon}_{\G}=0$ if $d>D$.
If $\G$ is distance-regular,  then using Theorem \ref{drdg3}, $\G$ is a connected normal digraph and $\tilde{\varepsilon}_{\G}=\varepsilon_{\G}=\varepsilon_{d}$. Now let $\G$ be a connected normal digraph and $\tilde{\varepsilon}_{\G}= \varepsilon_{d}$. We show that $\G$ is a distance-regular digraph. Clearly $D=d$ (since otherwise $\tilde{\varepsilon}_{\G}=0$ and $\varepsilon_{d}>0$, a contradiction) and
$\langle \tilde{\textbf{A}}_D,\tilde{\textbf{A}}_D \rangle=\frac{\langle \tilde{\textbf{A}}_D,P_D(\textbf{A}) \rangle^2}{\delta_D}$.
Hence $$\parallel \tilde{\textbf{A}}_D-Proj_{P_D(\textbf{A})}(\tilde{\textbf{A}}_D)\parallel=0.$$
Therefore
$$\tilde{\textbf{A}}_D=Proj_{P_D(\textbf{A})}(\tilde{\textbf{A}}_D)=\frac{\langle \tilde{\textbf{A}}_D,P_D(\textbf{A}) \rangle}{\delta_D}P_D(\textbf{A}).$$
Using an argument similar to the proof of Theorem \ref{drdg3}, we have $\tilde{\textbf{A}}_k=\frac{\langle H(\textbf{A}),P_{k}(\textbf{A}) \rangle}{\delta_{k}}P_{k}(\textbf{A})$ for each $0\leq k\leq D$ and hence using Theorem \ref{weakly distance-regular2}, $\G$ is a distance-regular digraph.\epf


\section{Connected normal digraphs with finite odd-girth}

In this section, using Theorem \ref{drdg9}, we demonstrate that a connected normal digraph $\G$  is either a bipartite digraph or a generalized odd graph or it has a finite odd-girth $g_o(\G)\leq \min\{2d-1,2D+1\}$.

\begin{lemma}\label{dog0} Let $\G$ be a connected digraph with finite odd-girth $g_o(\G)$. Then $g_o(\G)\leq 2D+1$.
\end{lemma}

\pf Assume to the contrary that $g_o(\G)> 2D+1$.
\begin{claim} The length of each walk of length at most $D+2$ from a vertex $u$ to a vertex $v$ is even if and only if $\partial(u,v)$ is even.
\end{claim}
To give a proof for the claim first assume that $W_1$ is a walk of odd length at most $D+2$ from $u$ to $v$ and $P$ of even length is the shortest path from $u$ to $v$. Now assume that $Q$ is the shortest path from $v$ to $u$, the existence of such a path is guaranteed by the fact that $\G$ is connected. For even $\partial(v,u)$, $W_1\cup Q$  and for odd $\partial(v,u)$, $P\cup Q$ is an odd closed walk of length at most $2D+1$ and so there is an odd cycle of length at most $2D+1$, a contradiction. A similar argument is used when $W_1$ is of even length and $P$ is of odd length.

Now let $C=v_1v_2\ldots v_{g_o(\G)}v_1$ be the shortest odd cycle and let $P$ be the shortest path from $v_1$ to $v_{D+2}$. Clearly $W=v_1Pv_{D+2}v_{D+3}Cv_{g_o(\G)}v_1$
is a closed walk of odd length at most $g_o(\G)-2$, a contradiction.
\epf

\begin{lemma}\label{dog1} Let $\G$ be a connected normal digraph with finite odd-girth $g_o(\G)\geq 2d+1$. Then $\G$ is distance-regular and $g_o(\G)=2d+1$.
\end{lemma}

\pf Since $D\leq \hat{D}=d$ ($\G$ is normal) and $2d+1\leq g_o(\G)\leq 2D+1$ (by Lemma \ref{dog0}), we have $D=d$ and $g_o(\G)=2d+1$. Now by Theorem \ref{drdg9}, it suffices to show that $\tilde{\varepsilon}_{\G}=\varepsilon_{d}$ (that is, the weighted excess is equal to the spectral excess). To do this it suffices to prove that $\langle \tilde{\textbf{A}}_D,\tilde{\textbf{A}}_D \rangle=\frac{\langle \tilde{\textbf{A}}_D,P_D(\textbf{A}) \rangle^2}{\delta_D}.$

An argument similar to the proof of Lemma \ref{dog0} is used to see that the length of each walk of length at most $D$ from a vertex $u$ to a vertex $v$ is even if and only if $\partial(u,v)$ is even. It follows that for every $0\leq i\leq j\leq D$, if $\langle A^i,A^j\rangle$ is non-zero, then both $i$ and $j$ are even or odd. It is easy to show by induction on $i$ that $P_i$ is an even or odd polynomial depending on whether $i$ is even or odd, for all $i\leq D$.

\begin{claim} \label{c2} For every $1\leq i\leq [D/2]$ there are no two vertices $u$ and $v$ with two $uv$-walks $W_1$ and $W_2$ of lengths $D+1$ and $D-2i$.
\end{claim}
Assume to the contrary that for some $1\leq i\leq [D/2]$ there are two vertices $u$ and $v$ with two $uv$-walks $W_1=v_1(=u)v_2\ldots v_{D+2}(=v)$ and $W_2$ of lengths $D+1$ and $D-2i$, respectively. Now let $P=u_1(=v)u_2\ldots u_{l+1}(=v_2)$ be the shortest path from $v$ to $v_2$. Clearly $l+D$ is even, since otherwise $vPv_2v_3W_1v_{D+2}$ is a closed walk of odd length at most $2D-1$, a contradiction. Since $u$ is adjacent to $v_2$ and $W_3=uW_2vPv_2$ is a walk of even length between $u$ and $v_2$ we have $l>2i\geq 2$ (since otherwise either $uW_3v_2Qu$ or $uv_2Qu$ is a closed walk of odd length at most $2D-1$ where $Q$ is the shortest path from $v_2$ to $u$, a contradiction). Now let
$R$ of length $r$ be the shortest path from $v_2$ to $u_2$. Clearly $r+l-1$ is even and so $r<D$. Hence either $uv_2Ru_2Su$ or $uW_2vu_2Su$  is a closed walk of odd length at most $2D-1$ where $S$ is the shortest path from $u_2$ to $u$, a contradiction. Hence for every $1\leq i\leq [D/2]$ we have $\langle \textbf{A}^{D+1},\textbf{A}^{D-2i} \rangle=0$.

Consider the vertices $u$ and $v$ with $\partial (u,v)=D$. Since $H(x)$ has degree $D$, we have $(\textbf{A}^{D})_{uv}=\zeta (H(\textbf{A}))_{uv}$ and since $(\textbf{A}-\lambda_0I)H(\textbf{A})=0$, we have $(\textbf{A}^{D+1})_{uv}=\gamma(H(\textbf{A}))_{uv}$, where $\zeta=\frac{\pi_0}{n}$ and $\gamma=\frac{\pi_0}{n}\Sigma_{i=0}^{D}\lambda_i$ are independent to the choice of $u$ and $v$. Assume that  $\zeta_D$ is the leading coefficient of $P_D$. Using Claim \ref{c2} we have:

$$\langle xP_D,P_D \rangle=\zeta^2_D\langle \textbf{A}^{D+1},\textbf{A}^{D} \rangle=\frac{\zeta^2_D\sum_{u\in V(\G)}(\textbf{A}^{D}(\textbf{A}^t)^{D+1})_{uu}}{n}=$$$$\frac{\zeta^2_D\gamma\zeta\sum_{u\in V(\G)}\sum_{v\in \G^{+}_D(u)}(H(\textbf{A}))^2_{uv}}{n}=\zeta^2_D\gamma\zeta\tilde{\delta}_D.$$
On the other hand, $$(H(\textbf{A}))_{uv}=\frac{\langle H,P_D \rangle}{\delta_D}(P_D(\textbf{A}))_{uv}=\frac{\langle H,P_D \rangle}{\langle xP_D,P_D \rangle}(\textbf{A}P_D(\textbf{A}))_{uv}=\frac{\zeta_D\gamma(H(\textbf{A}))_{uv}\langle H,P_D \rangle}{\zeta^2_D\gamma\zeta\tilde{\delta}_D}.$$
Therefore, $\langle H,P_D \rangle=\zeta_D\zeta\tilde{\delta}_D.$
It is easy to see that  $\langle \tilde{\textbf{A}}_D,P_D(\textbf{A}) \rangle=\zeta_D\zeta\tilde{\delta}_D.$ On the other hand,
$$H(\textbf{A})=\tilde{\textbf{A}}_0+\tilde{\textbf{A}}_1+\cdots+\tilde{\textbf{A}}_D=\sum_{i=0}^{D}\frac{\langle H,P_i \rangle}{\delta_i}P_i(\textbf{A}).$$
Therefore
$\langle \tilde{\textbf{A}}_D,\tilde{\textbf{A}}_D \rangle=\frac{\langle H,P_D \rangle\langle \tilde{\textbf{A}}_D,P_D(\textbf{A}) \rangle}{\delta_D}.$
Hence $\tilde{\delta}_D=\frac{\delta_D}{(\zeta_D\zeta)^2}$ and so
$$\langle \tilde{\textbf{A}}_D,\tilde{\textbf{A}}_D \rangle=\frac{\langle \tilde{\textbf{A}}_D,P_D(\textbf{A}) \rangle^2}{\delta_D}=\tilde{\delta}_D.$$
\epf


The following facts are from \cite{D}.

\begin{theorem}\label{drgg}{\cite{D}} Let $\G$ be a distance-regular digraph with girth $g$. Then
\begin{itemize}\item[\rm(i)]
If $0<\partial(u,v)<g$, then $\partial(v,u)=g-\partial(u,v)$;
\item[\rm(ii)]  $D=g$ or $D=g-1$;
\item[\rm(iii)]  If $\G$ satisfies $D= g- 1$ and $\textbf{A}$, is its adjacency matrix, then
$\textbf{A} \otimes J_m$, for any $m\geq 2$ is the adjacency matrix of a distance-regular digraph
with $D = g$, and conversely, any such digraph satisfying $D = g$ comes about
this way from one with $D= g - 1$.
\end{itemize}
\end{theorem}

\begin{theorem}\label{dog2} Let $\G$ be a connected normal digraph. Then one of the following holds:

\begin{itemize}
\item[\rm(i)] $\G$ is a bipartite digraph;
\item[\rm(ii)] $\G$ is a generalized odd graph;
\item[\rm(iii)]$g_o(\G)\leq \min\{2d-1,2D+1\}$.
\end{itemize}

\end{theorem}

\pf If $\G$ has no odd cycle, then $\G$ is bipartite; see \cite{BJG} . Now let $\G$ have finite odd-girth $g_o(\G)$. By Lemma \ref{dog0}, $g_o(\G)\leq 2D+1$. Now let $g_o(\G)\geq 2d+1$.
Using Lemma \ref{dog1}, $\G$ is distance-regular and $g_o(\G)=2d+1=2D+1$. If $g=2$, then $\G$ is a generalized odd graph. Now let $g\geq 3$. We will show that there is no such digraph. If $\G$ is long; that is, $D=g$, then a distance-regular digraph $\overline{\G}$ obtained by identifying all antipodal vertices of $\G$ has diameter $D-1$ (in fact, we have $\textbf{A}_{\G}=\textbf{A}_{\overline{\G}} \otimes J_m$ for some $m\geq 2$ and so  by Theorem \ref{drgg} $\overline{\G}$ has diameter $g-1=D-1$) and odd-girth $2D+1$, which is impossible by Lemma \ref{dog1}. Therefore, $g=D+1$ is even and so $D$ is odd. Let $C=v_1v_2\ldots v_{2D+1}v_1$ be the smallest odd cycle of $\G$. Since $v_{g+1}v_{g+2}\ldots v_1$ is a path of length $D$, an argument similar to the proof of the claim in Lemma \ref{dog0} implies that $\partial(v_{g+1},v_1)$ is odd. First let $\partial(v_1,v_{g+1})>1$. Using Theorem \ref{drgg} $\partial(v_1,v_{g+1})+\partial(v_{g+1},v_1)=g$ and so $\partial(v_1,v_{g+1})$ is odd. Hence $\partial(v_{g+1},v_1)<g-1$ and $v_{g+1}Pv_1v_2\ldots v_{g+1}$ is an odd closed walk of length at most $2g-2$ where $P$ is the shortest path from $v_{g+1}$ to $v_1$. Therefore, $g_o(\G)\leq 2g-2=2D$, a contradiction. So $\partial(v_1,v_{g+1})=1$. With the same argument we have $\partial(v_i,v_{g+i})=1$ for every $1\leq i\leq 2D+1$ where we use mod $2D+1$ arithmetic. Now $v_1v_{g+1}v_{g+2}\ldots v_{2D+1}v_{g}v_{1}$ is an odd cycle of length $D+2\leq 2D+1.$ Hence $g=D+1=2$, which is impossible.
\epf

The equality in Part (iii) of Theorem \ref{dog2} is attained for strongly regular diraphs with odd-girth three (in fact, for these digraphs we have $g_o(\G)=2d-1<2D+1$). On the other hand, among the triangle-free connected regular graphs with diameter two there are many graphs that are not strongly regular. These graphs have at least four distinct eigenvalues and so those which are not bipartite have odd-girth five (in fact, for these graphs we have $g_o(\G)=2D+1\leq 2d-1$). Hence the bound given in (iii) is the best possible.\\

As a corollary of Theorem \ref{dog2} we conclude that every connected normal digraph $\G$ with finite odd-girth at least $2d+1$ is a generalized odd graph, which generalizes the main result in \cite{vH} and also a result in \cite{LW}.

\begin{corollary}\label{dog3} Let $\G$ be a connected normal digraph with finite odd-girth $g_o(\G)\geq 2d+1$. Then $\G$ is  a generalized odd graph.
\end{corollary}

The odd-girth $g_o(\G)$ is the minimum odd number $k$ for which $tr(\textbf{A}^k)=\sum_{i=0}^{d}m_i\lambda^k_i$ is non-zero. Therefore, the odd-girth of a digraph is determined according to the spectrum and so by Corollary \ref{dog3} we have the following result, which implies a known result in \cite{HLiu} on the spectral characterization of generalized odd graphs.

\begin{corollary} Any connected normal digraph cospectral with a generalized odd graph, is a generalized odd graph.
\end{corollary}

\end{document}